\documentstyle[11pt,amscd,pstricks, pst-node,righttag]{amsart}
\psset{nodesep=5pt,arrows=->,linewidth=.4pt}

\allowdisplaybreaks

\newcommand{\id}{\operatorname{id}}
\newcommand{\Uq}{U_q(g)}
\newcommand{\Uz}{U_\zeta(g)}
\newcommand{\Gq}{G_q^{(\pi)} }
\newcommand{\Quea}{Quantized universal enveloping algebra }
\newcommand{\quea}{quantized universal enveloping algebra }

\begin{document}
\title[Quantum groups and Cuntz algebra]{Quantum Group Duality and the Cuntz Algebra}
\author{A. L. Carey}
\address[A. L. Carey]{Department of Pure Mathematics, University of Adelaide,
Adelaide, Australia}
\author{A. Paolucci}
\address[A. Paolucci]{Department of Mathematics, University of Leeds,
Leeds, LS2 9JT, UK}
\author{R. B. Zhang}
\address[R. B. Zhang]{Department of Pure Mathematics\\ University of Adelaide\\
Adelaide, Australia}
\date{}

\begin{abstract}
The Cuntz algebra carries in a natural way the structure of a module
algebra over the quantized universal enveloping algebra $U_q(g)$,
and the structure of a co-module algebra over the quantum group
$G_q$ associated with $U_q(g)$.  These two algebraic structures
are dual to each other via the duality between $G_q$ and $U_q(g)$.
\end{abstract}

\maketitle

\section{Introduction}

This paper is motivated by earlier work of one of us \cite{P}
on co-actions of (Woronowicz' \cite{W}) quantum groups on Cuntz algebras
\cite{Cuntz} and by studies by Doplicher and Roberts \cite{DR}
of group actions on Cuntz algebras.

The general theory of Hopf algebras \cite{M}
suggests that we should seek actions of the Drinfeld-Jimbo
\cite{Drinfeld,J} deformations of universal enveloping algebras
of semi-simple Lie algebras dual to the co-actions (found in \cite{P})
of the Woronowicz quantum groups.
This is achieved in this paper in a natural way.
An interesting feature is the appearance of braid groups
as the generators of the fixed and co-fixed points.
By analogy with \cite{DR} this paper provides
interesting parallels with the work
on braid group statistics in two dimensional
quantum field theory \cite{F}.

To summarize our results we collect them here
in the introduction into two theorems.
To make the proofs easily
accessible we include some expository material in  Section 2
where we introduce
notation and basic facts. We denote by
$\Uq$ the quantized universal enveloping algebra of a simple
Lie algebra $g$. When $0<q<1$ the algebraic dual of $\Uq$
contains the Woronowicz quantum group  $\overline{G^{(\pi)}_q}$.
(The bar indicates closure in Woronowicz' $C^*$-norm of a dense Hopf
subalgebra $G^{(\pi)}_q$ introduced in subsection 2.2).
In Section 3 we describe realizations of the braid group
in the Cuntz algebra
for generic $q$ (that is, not a root of unity).

When $0<q<1$ a co-action of $\overline{G^{(\pi)}_q}$ on the Cuntz algebra
$O_d$ on $d$ generators (where $d$ is the dimension of a
representation of $g$ chosen so that its tensor powers contain
all the irreducible representations of $g$ as subrepresentations)
was discovered in \cite{P}. It was also shown that there is a homomorphism
$\alpha$ of the braid group $B_n$ on $n$ generators into $O_d$
such that the co-fixed points under the $\overline{G^{(\pi)}_q}$ co-action
were generated by $\alpha(B_n)$. In section 3 we review this
embedding of the braid groups in $O_d$.

The main results proved in this paper in Sections 4 and 5
may be stated as follows.

\noindent
{\bf Theorem 1}: {\em For generic $q$
there is an action of $\Uq$ on the Cuntz algebra
by non-unital non-star densely defined endomorphisms with the following
properties:\\
(i) The domain of the endomorphisms is the dense subalgebra $O_d^0$
of $O_d$ given by polynomials in the generators.\\
(ii) There is a co-action
of $G^{(\pi)}_q$ on  $O_d^0$ dual to
this $\Uq$ action. This co-action coincides for $0<q<1$ with
the co-action of $\overline{G^{(\pi)}_q}$ discovered in
\cite{P};  \\
(iii) The braid
group elements in $O_d^0$ are both fixed and co-fixed.\\
(iv) The  fixed and co-fixed elements of $O_d^0$ coincide.}

We relate this to the results in \cite{P} at the end of Section 5
in the following way.

\noindent
{\bf Theorem 2}: {\em (i) There is a distinguished element of $O_d^0$ which we
call the rank $d$ $q$-antisymmetric tensor which is fixed and co-fixed
by the action and co-actions of theorem 1.\\
(ii) When $0<q<1$ the  braid group elements together with
the rank $d$ $q$-antisymmetric tensor
algebraically generate the fixed points of the $\Uq$
action and generate, in the topological sense,
 the co-fixed subalgebra of $O_d$ under the  $\overline{G^{(\pi)}_q}$ action.}

Section  4 contains the definition of the $\Uq$ action, and also
discusses the braids as fixed points. Section 5 gives the dual
co-action of $G^{(\pi)}_q$ and completes the proof of both theorems.
 A natural question to ask is
what happens if one dualizes the co-action and section 6 is devoted
to showing that we recover the given action of $\Uq$.
There is a slight technical point in that the Hopf dual of $G^{(\pi)}_q$
is not obviously $\Uq$. There is a possibly non-trivial ideal
which one needs to factor out. However in the Cuntz realization this
ideal acts trivially so that we have duality working fully.

In order to understand the observations of this paper
at a more fundamental level we conclude our discussion in
the final section by developing the action and co-action within
the framework of braided tensor categories (cf \cite{P}).
The viewpoint in this paper is algebraic in that we have not attempted
to analyse the role played by the $C^*$-algebra topology on the Cuntz
algebra from the viewpoint of the $\Uq$ action.
 It may be that there is an interesting connection between
this topology, that on $G^{(\pi)}_q$ described in \cite{W},
and some, as yet undetermined, topology on $\Uq$.

\section{Quantized universal enveloping algebra and dual quantum group}

We review the definitions and some properties of the quantized
universal enveloping algebras and their dual quantum groups in this section.

\subsection{\Quea \  $\Uq$}
Let $g$ be any finite dimensional simple Lie algebra over the complex
field $\Bbb C$. Denote by $\Phi^+$ the set of the positive roots of
$g$ relative to a base  $\Pi=\{\alpha_1, ..., \alpha_r\}$, 
where $r$ is the rank of $g$.
Define $E=\bigoplus_{i=1}^r{\Bbb R}\alpha_i$.    Let $(\ , \ ): E\times E
\rightarrow {\Bbb R}$ be an inner product of $E$  such that the Cartan
matrix $A$ of $g$ is given by
\begin{align*}
A=\left( a_{i j}\right)_{i j=1}^r, & \quad 
a_{i j}={ {2(\alpha_i, \alpha_j)}
\over{(\alpha_i, \alpha_i)}}.
\end{align*}
The Jimbo version  \cite{J} of the   \quea  $\Uq$ is defined to be the
unital associative algebra over $\Bbb C$,  generated by $\{ h_i, e_i, f_i\ | \
i=1, ..., r\}$  with the following relations
\begin{align}
k_i k_j = k_j k_i, & k_i k_i^{-1} =1, \nonumber\\
k_i e_j k_i^{-1} = q_i^2 e_j, &
k_i f_j k_i^{-1} = q_i^{-2} f_j, \nonumber
\end{align}
\begin{align}
{[}e_i, \ f_j]=\delta_{i j}{{ k_i - k_i^{-1}}
\over{q_i -q_i^{-1}}}, &\nonumber\\
\sum_{t=0}^{1-a_{i j}} (-1)^t
\left[ \begin{array}{r r r}
        1&-&a_{i j} \\
           &t&
        \end{array}\right]_{q_i}
(e_i)^t e_j (e_i)^{1-a_{i j}-t}&= 0 , \ \ \ \ i\ne j,  \nonumber\\
\sum_{t=0}^{1-a_{i j}} (-1)^t
\left[ \begin{array}{r r r}
        1&-&a_{i j} \\
           &t&
        \end{array}\right]_{q_i}
(f_i)^t f_j (f_i)^{1-a_{i j}-t}&= 0 , \ \ \ \ i\ne j, \label{quea}
\end{align}
where  $q$ is a complex parameter, which is assumed to be non-zero,
and is not a root of unity.
Also, $\left[ \begin{array}{r}
       s\\ t
        \end{array}\right]_q $
is the Gauss polynomial, and  $q_i=q^{(\alpha_i, \ \alpha_i)/2}$.

The algebra $\Uq$ has in addition the structure of a Hopf algebra.
We take the following co-multiplication
\begin{align*}
\Delta(k_i^{\pm 1})&=k_i^{\pm 1}\otimes k_i^{\pm 1},\\
\Delta(e_i)&=e_i\otimes k_i + 1\otimes e_i, \\
\Delta(f_i)&=f_i\otimes 1 + k_i^{-1}\otimes f_i.
\end{align*}
The co-unit $\epsilon: \Uq \rightarrow {\Bbb C}$
 and antipode $\gamma: \Uq \rightarrow \Uq$
 are respectively given by
\begin{align*}
 \epsilon(e_i) =\epsilon(f_i)=0,\qquad  &\epsilon(k_i^{\pm 1})= \epsilon(1)=1,
\end{align*}
\begin{align*}
\gamma(e_i)&=-e_i k_i^{-1}, \\
\gamma(f_i)&=-k_i f_i, \\
\gamma(k_i^{\pm 1})&=k_i^{\mp 1}.
\end{align*}

Later in the paper we will need the notion of a `universal 
$R$-matrix' in the sense used in connection with
 the Yang-Baxter equations. This $R$-matrix
does not live naturally in the Jimbo picture.
For this we need the Drinfeld version \cite{D}.
If we set
\begin{align*}
q& = \exp(\zeta),
\end{align*}
and regard $\zeta$ as a formal indeterminate,
then the Drinfeld version \cite{Drinfeld} of the \quea is an  associative
algebra over ${\Bbb C}[[\zeta]]$ completed with respect to the
$\zeta$-adic topology for ${\Bbb C}[[\zeta]]$.
It is generated by $\{e_i, \ f_i, \ h_i, \ i=1, 2, ..., r\}$ with
\begin{align*}
k_i^{\pm 1}&=q^{\pm h_i},
\end{align*}
subject to the same relations (\ref{quea}).

We use the notation $U_\zeta(g)$
to denote this algebra which in
 the terminology of Drinfeld is a quasi-triangular Hopf
algebra. This means it admits an invertible
$R\in U_\zeta(g){\hat\otimes} U_\zeta(g)$ ( ${\hat\otimes}$ 
represents tensor product completed with respect to the
$\zeta$ - adic topology ), called the universal
$R$-matrix, which satisfies the following defining relations
\begin{align*}
R\Delta(a)&=\Delta^\prime(a)R,\qquad \forall a\in \Uz,\\
(\Delta\otimes\id)R&=R_{13}R_{23},\qquad (\id\otimes\Delta)R=R_{13}R_{12}.
\end{align*}
Further general properties of $R$ are
\begin{align*}
(\gamma\otimes\id)R&=(\id\otimes\gamma)R=R^{-1}\\
(\epsilon\otimes\id)R&=(\id\otimes\epsilon)R=1,\\
%\noalign{and}
R_{12}R_{13}R_{23}&=R_{23}R_{13}R_{12},
\end{align*}
where the last equation is the celebrated quantum Yang-Baxter
equation, which is a direct consequence of the defining relations of
$R$.

The universal $R$-matrix is of the form
\[
R=q^{\sum_{i,j}(B^{-1})_{ij}h_i\otimes h_j}[1\otimes 1+\sum_s C_{(s)}E^{(s)}
  \otimes F^{(s)}],
\]
where
\begin{align*}
&B=((\alpha_i, \alpha_j))\\
&E^{(s)}\mbox{ are combinations of products of } e_i\mbox{'s},\\
&F^{(s)}\mbox{ are combinations of products of } f_i\mbox{'s},\\
&C_{(s)}\mbox{ are scalars in } {\Bbb C}[[\zeta]].
\end{align*}
Moreover,
\begin{align} { {R -1\otimes 1}\over{\zeta}}
&= \sum_{i,j}(B^{-1})_{ij}h_i\otimes h_j +
\sum_{\alpha\in\Phi^+} e_\alpha\otimes f_\alpha  + o(\zeta),
\label{T}
\end{align}
where $e_\alpha$ and $f_\alpha$ are the quantum analogs of the Cartan-Weyl
generators.

Given any nontrivial irreducible representation $\pi$ of $\Uz$,
we set
\begin{align}
T^{(+)}=(\pi \otimes id ) R^T, & \qquad T^{(-)}=(\pi \otimes id ) R^{-1}.
\end{align}

\noindent{\bf Lemma 1}: {\it The matrix elements of $T^{(+)}$
and $T^{(-)}$ generate 
the entire Drinfeld algebra $U_\zeta(g)$ topologically.}

\noindent {\bf Proof}.
Recall that the Drinfeld algebra is  a deformation of
the universal enveloping algebra $U(g)$ of the simple Lie
algebra $g$ in the sense of Gerstenhaber \cite{Gerstenhaber}.
Therefore
\begin{align*}
U(g)&= \Uz/\zeta \Uz.
\end{align*}
Denote by $U(g)[[\zeta]]$ the universal enveloping algebra of $g$ over
${\Bbb C}[[\zeta]]$ completed with respect to the $\zeta$-adic topology.
The rigidity of $U(g)$ leads to the conclusion that
$U(g)[[\zeta]]$ and $U_\zeta(g)$ are isomorphic
as associative algebras and hence they have the same
representation theory. Given
 these facts we can easily deduce that in order to prove the Lemma,
it suffices to show that the matrix elements of
\begin{align*}
{{1}\over{\zeta}}\left( T^{(\pm)} - I\otimes 1\right)/&\zeta U_\zeta(g)
\end{align*}
generate $U(g)$. Here $\zeta U_\zeta(g)$ denotes the ideal in
 $U_\zeta(g)$ generated by $\zeta$.
Using the well known fact that the trace over any finite dimensional
irreducible representation of the simple Lie algebra $g$ defines
an invariant non-degenerate bilinear form, we obtain from (\ref{T})
\begin{align*}
tr_\pi\left\{ {{1}\over{\zeta}}\left(T^{(\pm)} - I\otimes 1\right)
     [\pi(h_i)\otimes 1]\right\}& = \pm c_\pi {\hat h}_i + o(\zeta), \\
tr_\pi\left\{ {{1}\over{\zeta}}\left(T^{(+)} - I\otimes 1\right)
     [\pi(e_i)\otimes 1]\right\}& =  c_\pi e_i + o(\zeta), \\
tr_\pi\left\{ {{1}\over{\zeta}}\left(T^{(-)} - I\otimes 1\right)
     [\pi(f_i)\otimes 1]\right\}& = - c_\pi f_i + o(\zeta),
\end{align*}
where $c_\pi\in {\Bbb C}$ is a nonvanishing constant, and the ${\hat h}_i$
are independent linear  combinations of the $h_j$ such that
\begin{align*}
{[}{\hat h}_i,  e_j]&=a_{i j} e_j.
\end{align*}
Quotienting by  $\zeta U_\zeta(g)$ is equivalent to taking
 $\zeta\rightarrow 0$ in the preceeding equations
whose right hand sides 
then clearly generate $U(g)$, thus completing the proof of the
Lemma.

Observe also the following important fact:  $T^{(\pm)}$ can be expressed
solely in terms of $e_i, \ f_i, \ k_i^{\pm 1}$,
and  $k_\pi^{\pm 1} = q^{\pm h_\mu}$, with $\mu$ being the
highest weight of $\pi$,
and $h_\mu$ a linear combination of the $h_i$ such that
$[h_\mu, e_i] = (\mu, \alpha_i) e_i$.  For some Lie algebras it can happen
that $\mu$ is not in the root lattice of $g$. In that case,
$k_\pi^{\pm 1}$ can not be expressed as products of integer powers of
the $k_i^{\pm 1}$.

Let us now  specialize  $q$ to a complex parameter.
The matrix elements of $T^{(\pm)}$ generate, algebraically, 
an associative algebra ${\tilde U}_q(g)$ over the complex field $\Bbb C$. 
If $\mu$ belongs to the root lattice of $g$,  
this algebra coincides with the Jimbo \quea, 
otherwise, it contains $\Uq$ as a subalgebra.
The action of the Jimbo algebra we define
in section 4 extends to the  algebra ${\tilde U}_q(g)$ but not uniquely
in the latter case.
The extensions differ however only by the fact that one needs to
choose a particular (complex) root of $q$, and are related  
to one another by an action of the appropriate group of
 roots of unity as automorphisms.
At generic $q$, this difference
between the two algebras $\Uq$ and ${\tilde U}_q(g)$ 
is thus not at all important.
Finally we remark
that if we further assume that $q$ is real, then $U_q(g)$ has the structure
of a Hopf *-algebra.

\subsection{ The quantum group  $G^{(\pi)}_q$}
We now move on to set up the notation and properties
of the dual to $\Uq$, the Woronowicz quantum group.
The finite dual $(\Uq)^0$ of $\Uq$ has a natural Hopf algebra
structure, with the multiplication $m_0$, co-multiplication $\Delta_0$,
unit ${\bf 1}_0$, co-unit $\epsilon_0$, and  antipode $\gamma_0$
defined in the standard fashion \cite{M}.
  We consider a subalgebra $G^{(\pi)}_q$
of  $(\Uq)^0$ defined in the following way.
Let $\pi$ be a finite dimensional non-trivial representation of
$\Uq$, which can be assumed to be irreducible without losing generality.
Set $d= dim_{\Bbb C}\pi$.  Consider the matrix
\begin{align*}
         U= \left( u_{i j}\right)_{i, j =1}^d,  &\quad 
 u_{i j}\in (\Uq)^0,
\end{align*}
defined by
\begin{align*}
   \langle u_{ij}, a \rangle &= \pi(a)_{ij},  \ \ \ \ \forall a\in\Uq.
\end{align*}

\noindent{\bf Definition}: We define
$G^{(\pi)}_q$ as the associative subalgebra of $(\Uq)^0$
generated by the matrix elements of $U$, with the
multiplication  defined by
\begin{align*}
%\intertext{and}
   \langle u_{ij}u_{kl}, a\rangle &= \sum_{(a)} \pi_{ij}(a_{(1)})\pi_{kl}
      (a_{(2)}),  \ \ \ \ \forall a\in\Uq.
\end{align*}

For $0<q<1$ in \cite{W}
 this algebra is completed in an appropriate $C^*$-algebra
norm. For the most part we will not need this topology here and will work with
the uncompleted algebra. When necessary we denote the $C^*$ completion
by a bar.
Clearly the unit of $\Gq$ coincides with the co-unit $\epsilon$ of $\Uq$.
Set
\[
   R^{(\pi)} = (\pi\otimes\pi) R,
\]
where $R$ is the universal $R$-matrix of $\Uq$.

\noindent{\bf Lemma 2}:    {\it 
$U$ satisfies the following quadratic relation}
\begin{equation}
   R^{(\pi)}_{12}U_1U_2 = U_2U_1R_{12}^{(\pi)}. \label{Gq}
\end{equation}

\noindent{\bf Proof}: We first
 note that the  left-hand side contracted with any
$a\in \Uq$ gives
\[
   R^{(\pi)}_{12}U_1U_2(a) = R^{(\pi)} (\pi\otimes\pi)\Delta(a);
\]
while the right hand side yields
\[
   U_2U_1R_{12}^{(\pi)}(a) = (\pi\otimes \pi)\Delta^\prime(a)R^{(\pi)}.
\]
Then equation (\ref{Gq}) immediately follows from the fact that
\[
  R^{(\pi)}(\pi\otimes \pi)\Delta(a) = (\pi\otimes\pi)\Delta^\prime(a)
  R^{(\pi)},\quad \forall a\in \Uq, 
\]
completing the proof.

It is also easy to see that $\Gq$ has the structure of a bi-algebra.
The co-multiplication is given by
\begin{equation}
   \Delta_0(u_{ij}) = \sum_k u_{ik}\otimes u_{kj},
\end{equation}
which follows from the equation
\begin{align*}
   \langle \Delta_0(u_{ij}), a\otimes b \rangle &= \pi_{ij}(ab), \ \ \ \
\forall a, b \in\Uq,
\end{align*}
while the co-unit is $1_{\Uq}$
\[
  u_{ij}(1_{\Uq}) = \delta_{ij}.
\]
Furthermore,  $\Gq$  admits an antipode $\gamma_0$.
Let $\pi^\dagger$ be the irreducible representation of $U_q(g)$ dual to
$\pi$.  If $V$ is the module furnishing the representation $\pi$, 
we take $V^*$ to be the dual vector space on which
 the dual representation is defined by
$(\pi^\dagger(a).v^* ) ( v') = v^*( \pi( \gamma(a) ).v' )$, 
where  $\gamma$ is the antipode for $\Uq$.  
Then the antipode for  $\Gq$ is given by
\begin{align}
   \gamma_0(U)& = U^{-1},
\end{align}
where $U^{-1}$ is another $d\otimes d$ matrix, the entries of which
belong to $(U_q(g))^0$ and satisfy
\begin{align*}
   (U^{-1})_{i j}( a ) &= \pi^\dagger( a )_{i j}, \quad   \forall a\in U_q(g).
\end{align*}
For this to make sense, we need to show that

\noindent{\bf Lemma 3}:
{\it The elements of $U^{-1}$ belong to $\Gq$.}

\noindent{\bf Proof}. The Lemma is equivalent to the statement that
some repeated tensor product of $\pi$
 (with respect to the
co-multiplication
$\Delta$) contains the dual representation $\pi^\dagger$ as an
irreducible component.  For this to be true, it suffices to show
that a one dimensional representation can arise from
nontrivial tensor products of $\pi$. Let $V$ be the $U_q(g)$ module
which furnishes the representation $\pi$. We claim that there exists
a nonvanishing $\Lambda\in V^{\otimes d}$ which generates a one
dimensional representation of $U_q(g)$.
As the representation theory
of $U_q(g)$ at generic $q$ is the same as that of $U(g)$, let us
first examine the classical situation.
In the $q\rightarrow 1$ limit, the $U_q(g)$ action on $V$ yields
a $U(g)$ action.  The totally antisymmetric rank $d$ tensor of
$V$ is one dimensional, and we denote its basis element by $\Lambda^{(0)}$.
It is clearly true that
\begin{align*}
x\cdot \Lambda^{(0)} &= 0, \quad \forall x\in g.
\end{align*}
By calling upon the Lusztig-Rosso theorem, \cite{L,R} we conclude that
there exists a
nonvanishing  $\Lambda\in V^{\otimes d}$, which reduces to
$\Lambda^{(0)}$ in the $q\rightarrow 1$ limit,  such that
\begin{align*}
a\circ \Lambda &= \epsilon(a), \quad \forall a\in U_q(g).
\end{align*}
This completes the proof of the Lemma.

\noindent{\bf Definition}. We will call $\Lambda$ the
rank $d(=\dim\pi)$ $q$-antisymmetric tensor of $V$.

\subsection{Modules and co-modules}

We begin with some generalities on co-actions for Hopf algebras.
Let $A$ be a Hopf algebra with
 co-multiplication $\Delta$, co-unit $\epsilon$ and
antipode $\gamma$. Let $V$ be a left  $A$-module, which is assumed to
be locally finite, i.e., it satisfies the folllowing properties:
corresponding to each $v\in V$, we can find a finite set of elements
$v_i\in V$, $i=1, 2, ..., N$, such that
\begin{align*}
v&=\sum_{i=1}^N c_i v_i, \quad c_i\in{\Bbb C}, \\
a\circ v_i& = \sum_{j=1}^N \psi_{j i}(a) v_j,  \quad \forall a\in A,
\end{align*}
where $\psi_{j i}(a)\in{\Bbb C}$.

Let $A^0$ be the finite  dual Hopf algebra with multiplication
$m_0$, unit $\epsilon$, co-unit $1_A$, co-multiplication $\Delta_0$
and antipode $\gamma_0$. The left $A$-module $V$ automatically carries
a right $A^0$ co-module structure
\[
\omega:V\rightarrow V\otimes A^0
\]
defined in the following way.
For any element $v\in V$, if we write
\[
 \omega(v)=\sum_{(v)}v_{(1)}\otimes v_{(2)}\qquad\qquad v_{(2)}\in A^0,
\]
then for all $a\in A$,
\[
 \omega(v)(a)=\sum_{(v)}v_{(1)}\ v_{(2)}(a)=a\circ v
\]

To be more explicit, we consider
\begin{align}
\omega(v_i)&=\sum_j v_j\otimes u_{ji},  \label{comodule}
\end{align}
where $u_{i j}$ are elements of $A^0$ which are uniquely
determined by the requirement that
\[
u_{i j}(a) = \psi(a)_{i j}, \ \ \ \ \ \forall a\in A.
\]
It follows from
\[
\{(id \otimes \Delta_0)\omega(v_i)\}(a\otimes b)
= \omega(v_i)(a b)
\]
that
\[ \Delta_0(u_{ij})=\sum_k u_{ik}\otimes u_{kj}.\]
Now
\[
(\id \otimes \Delta_0)\omega(v_i)=
\sum_{j,k}v_j\otimes u_{jk}\otimes u_{ki}.
\]
On the other hand,
\begin{align*}
(\omega\otimes \id)\omega(v_i)&=\sum_j \omega(v_i)\otimes u_{ji}
=\sum_{j,k} v_k\otimes u_{kj}\otimes u_{ji},
\end{align*}
hence
\[
(\id \otimes \Delta_0) \omega=( \omega\otimes \id) \omega.
\]
Denote by $\epsilon^*$ the co-unit of $A^0$, i.e., for any $a^*\in A^0$,
$\epsilon^*(a^*)=\langle a^*, 1_A\rangle$. Then
\[
(\id_V\otimes \epsilon^*) \omega(v_i)=\sum_j v_j\otimes \delta_{ji}=v_i,
\]
i.e.
\[
(\id_V\otimes \epsilon^*) \omega =\id_V.
\]
The above equations show that $\omega$ is indeed a right-co-module action
of $A^0$ on $V$.

\subsection{Fixed points and co-fixed points}
We continue the notation of the previous subsection.
As $V$ is a module over $A$, and at the same time, a co-module over $A^0$.
we introduce the notation
 $(V)^A$ for the fixed point set of $V$ under the action of $A$,
and $(V)_{A^0}$ for
the co-fixed point set of $V$ under the co-action of $A^0$.
Thus
\begin{align*}
(V)^A&=\{ v\in V \ | \  a\circ v=\epsilon(a)v,\  \forall a\in A\},  \\
(V)_{A^0} &=\{ u\in V \ | \ w(u)=u\otimes 1_{A^0}\}.
\end{align*}
Consider $v\in V$. If $v$ is co-fixed by $A^0$, i.e.,
\[
 \omega(v)=v\otimes \epsilon,
\]
then
\[
a\circ v=\omega(v)(a)= \epsilon(a)v,
\]
hence $v$ is also fixed by $A$.
On the other hand,  if $v\in V$ is fixed by $A$, then
\[
a\circ v=\epsilon(a) v=\sum_{(v)}v_{(1)}\ v_{(2)}(a),\quad \forall a\in A,
\]
i.e.
\[
 \omega(v)=v\otimes \epsilon,
\]
thus $v$ is also co-fixed. Therefore, $(V)^A=(V)_{A^0}$.

\subsection{The  $q$-determinant}
Returning now to the discussion 
begun in subsection 2.2 we will need the notion of $q$-determinant of
the matrix $U$ used to define $G_q^{(\pi)}$.
Let $V$ denote the representation space of $\pi$.
By the discussion of the preceeding two subsections we may consider
the co-action of $G_q^{(\pi)}$, denoted by $\omega$, on $V$. 
Now this co-action when applied to
the $q$-antisymmetric tensor $\Lambda$ yields
\[
\omega(\Lambda)=\Lambda\otimes\operatorname{det}_qU,
\]
where ${\det}_qU$ is some element of $\Gq$.

\noindent{\bf Definition}.   We  call ${\det}_qU$ the $q$-determinant
of $U$.

However, invoking the argument of the preceding
subsection, $\Lambda$ must be a co-fixed point of $\Gq$ as $\Lambda$
generates a trivial module of $U_q(g)$.  Hence
\[
\operatorname{det}_qU=\epsilon,
\]
i.e., $U$ has $q$-determinant $1_{G_q}(=\epsilon)$.

Note that in \cite{W} the quantum determinant is set equal to the identity
as an additional relation, while this relation is built into
 the definition of
$G_q^{(\pi)}$ used here.

\noindent{\bf Remarks}:
We need to comment on the representation $\pi$ and the dependence on it of
$G_q^{(\pi)}$.
  If $\pi$ can generate all the finite dimensional representations
of $U_q$ by repeated tensor
products, then $G_q^{(\pi)}$ will contain all $G_q^{(\pi^\prime)}$ as Hopf
subalgebras, where $\pi^\prime$ is any representation of $\Uq$.  The vector
representation
of $U_q(sl(n))$ has this property, as does the spinor representation
of $U_q(so(n))$.  A second point is that for a generic $q$, all finite
dimensional representations of $\Uq$ are
completely reducible, so we can simply take $\pi$ to be an irreducible
representation.  In general, repeated tensor products of $\pi$ can be
expressed as direct sums of a subclass of irreps of $U_q(g)$.
Components of these irreps form a basis of $G_q^{(\pi)}$, and this
statement may be regarded as an algebraic analogue of the
Peter-Weyl theorem in the quantum setting.
The multiplication rule of $G_q^{(\pi)}$ in terms of such a
basis is the Clebsch-Gordon decomposition of representations of $\Uq$.

\subsection{The dual Hopf algebra of $G_q^{(\pi)}$}

In this subsection we regard $G_q^{(\pi)}$ as a Hopf algebra in its
own right. It is defined to be the Hopf algebra generated by
the matrix elements of $U$, subject to the relations
(\ref{Gq}) together with the $q$-determinant condition
\[
\operatorname{det}_qU=1_{\Gq}.
\]
The co-multiplication, co-unit and antipode are as given before.
We want to investigate the finite dual of $\Gq$.

Consider the set of elements $\{l_{ij}^{(\pm)}\in(G_q^{(\pi)})^*|
i, j =1, ..., d\}$, which satisfy the following properties:
write
\[
L^{(\pm)}=\sum e_{ij}\otimes l_{ij}^{(\pm)}.
 \]
Then
\begin{align}
&\langle 1, U^{\otimes r}\rangle=I^{\otimes r}\nonumber \\
&\langle L^{(+)}, U_1U_2\ldots U_r\rangle=R_1^{(+)}R_2^{(+)}\ldots R_r^{(+)},
\nonumber \\
&\langle L^{(-)}, U_1U_2\ldots U_r\rangle=R_1^{(-)}R_2^{(-)}\ldots R_r^{(-)},
\label{definition}
\end{align}
where
\[
R^{(+)}=PR^{(\pi)}P,\qquad R^{(-)}=(R^{(\pi)})^{-1},
\]
and our notation is largely the same as that of \cite{Faddeev}.
By considering the co-multiplication of $G_q^{(\pi)}$, one can show
\cite{Faddeev} that $L^{(\pm)}$  satisfy the following relations
\begin{align}
R_{12}^{(+)}L_1^{(\pm)}L_2^{(\pm)}&=L_2^{(\pm)}L_1^{(\pm)}R_{12}^{(+)},
\nonumber \\
R_{12}^{(+)}L_1^{(+)}L_2^{(-)}&=L_2^{(-)}L_1^{(+)}R_{12}^{(+)}.
\end{align}

Let us denote by  $U_q^\prime$ the algebra generated by the
matrix elements of $L^{(\pm)}$.
A co-multiplication for $U_q^\prime$ is given by
\[
\Delta(L^{(\pm)})=L^{(\pm)}\otimes L^{(\pm)}\qquad (\mbox{i.e. }
\Delta(l_{ij}^{(\pm)})=\sum_k l_{ik}^{(\pm)}\otimes l_{kj}^{(\pm)}),
\]
and the corresponding co-unit and antipode are respectively
given by
\[
\epsilon(L^{(\pm)})=I, \qquad
\gamma(L^{(\pm)})=(L^{(\pm)})^{-1}.
\]

Define a linear map
\[  L^{(+)}\mapsto T^{(+)},
\ \ \  L^{(-)}\mapsto T^{(-)}. \]
This map extends in a unique way to a Hopf algebra homomorphism
$\phi:  U_q^\prime \rightarrow {\tilde U}_q(g)$, 
which is clearly surjective.  The kernel of $\phi$ is a Hopf ideal 
of $U_q^\prime$, and $U_q^\prime/Ker\phi\cong {\tilde U}_q(g)$ 
as Hopf algebras.

\section{Cuntz algebra realizations of the braid group}

\subsection{Representations of the braid group}
Consider $R^{(\pi)}$, where $\pi$ is a $d$-dimensional
representation of $\Uq$ and
let $P: V\otimes V\to V\otimes V$ be the flip:
$P(a\otimes b)= b\otimes a$. Define
$$
\sigma =PR^{(\pi)}.
$$
It follows from the defining relations of the universal $R$ matrix
that
\begin{align*}
&[\sigma , (\pi\otimes\pi )\Delta (a)]=0,\qquad \forall\ a\in \Uq, \\
  & (\sigma\otimes 1)(1\otimes\sigma)(\sigma\otimes 1) =
      (1\otimes\sigma)(\sigma\otimes1)(1\otimes\sigma)
\end{align*}
Note that $\sigma$ acts on $V\otimes V$,
while the above equation holds as endomorphisms of $V\otimes V\otimes V$.

Define
\[
  b_i = \underbrace{1\otimes\cdots\otimes 1}_{i-1}\otimes\sigma\otimes
 \underbrace{1\otimes\cdots\otimes 1}_{n-i},\qquad i=1,2,\ldots,n
\]
Then we have the braid relations
\begin{align}
   b_ib_{i+1}b_i &= b_{i+1}b_ib_{i+1} \nonumber \\
   b_ib_j &= b_jb_i,\qquad |i-j|\ge 2
\end{align}
Also,
\[
   \left[b_i, \pi^{\otimes(n+1)}\Delta^{(n)}(a)\right]=0,
    \qquad\forall\ a\in \Uq
\]

\subsection{Cuntz algebra realization of the braid group}
Let $O_d$ be the Cuntz algebra \cite {Cuntz} on $d$ generators.
This is the universal C$^*$-algebra generated by
$\{s_j | j=1,2,\ldots,d\}$ satisfying
$s_i^*s_j=\delta_{ij}$ and $\sum_j s_j s_j^* = 1$
where $1$ denotes the identity of $O_d$.

Let $\{e_{ij} | i,j = 1,2,\ldots,d\}$ be the matrix units for End $V$,
 which obey
\[
   e_{ij}e_{kl} = \delta_{jk}e_{il}.
\]
There exists a well known algebra homomorphism
\[   \eta:  \mbox{End} V \rightarrow  O_d, \]
defined by
$$   \eta(e_{ij})=s_is_j^*.
$$
This extends to an algebra homomorphism $\eta: \mbox{End} V^{\otimes m}
\rightarrow  O_d$ for each $m$, defined by
\begin{align*}
   &\eta\left(e_{i_1j_1}\otimes e_{i_2j_2}\otimes\cdots\otimes e_{i_mj_m}\right) \\
   &\qquad = s_{i_1}s_{i_2}\dots s_{i_{m-1}}s_{i_m}s^*_{j_m}s^*_{j_{m-1}}\dots
      s^*_{j_2}s^*_{j_1}.
\end{align*}

We notice that $\sigma\in\mbox{End} V$ can be expressed as
$$
   \sigma=\sum_{i,j, k, l=1}^d \sigma_{i j, k l}e_{ij}\otimes e_{kl},
$$
which,   under $\eta$, maps to
\begin{align}
   \theta = \eta(\sigma) = \sum \sigma_{ij,kl} s_i s_k s^*_l s^*_j.
\end{align}
More generally, we have
\begin{align}
   \theta_i=\eta(b_i)
 & = \sum_{\{l\}} s_{l_1}s_{l_2}\dots s_{l_{i-1}} \theta
     s^*_{l_{i-1}}\dots s^*_{l_2}s^*_{l_1},\nonumber\\
 & i = 1, 2, ..., n.
\end{align}
As elements of the Cuntz algebra, the $\theta_i$'s satisfy the braid
group relations.

\section{The Cuntz algebra as a module algebra over $U_q(g)$}

Introduce the notation:
\begin{align*}
  &H = \bigoplus_{i=1}^d {\Bbb{C}} s_i, \quad
   H^0={\Bbb C},  \quad
   H^r = H^{\otimes r}, \quad
   {\cal H} = \bigoplus^{\infty}_{r=0}H^r, 
\end{align*}
where  $\cal H$ is the  algebraic direct sum.
Let $H^*$ be the dual vector space of $H$,
\[
   H^* = \bigoplus_1^d {\Bbb{C}} s_i^*,
\]
with the pairing $H^*\otimes H\to\Bbb{C}$ given by
\[
   \langle s^*_i, s_j\rangle = \delta_{ij}.
\]
Set
\begin{align*}
 &  H^{*r}= H^{*\otimes r},\quad
 {\cal H}^* = \bigoplus_{r=0}^\infty H^{*r}.
\end{align*}
The pairing between $H^{*r},\ H^t$ is given by
\[
   \langle s^*_{j_r} s^*_{j_{r-1}}\dots s^*_{j_1}, s_{i_1} s_{i_2}\dots s_{i_t}
      \rangle = \delta_{rt}\delta_{i_1j_1}\delta_{i_2j_2}\dots\delta_{i_r j_r}.
\]
This definition is compatible with Cuntz multiplication when $r$ and $t$ are
equal. We now build up our action of $\Uq$ on the dense subalgebra
of $O_d$ consisting of polynomials in the generators in three steps
outlined in the following subsections. For convenience we introduce a 
notation for this subalgebra.

\noindent{\bf Definition}: Let $O_d^0$ denote the dense subalgebra
of $O_d$ consisting of polynomials in the generators.

\subsection{$\Uq$  action on $\cal{H}$ and $\cal{H}^*$}
We begin by defining the action of $U_q(g)$ on $\Bbb{C}$ by
\begin{align}
   a\circ c = \epsilon(a)c.
\end{align}
That is, we regard $\Bbb{C}$ as a trivial module over $\Uq$.
Now
let $\pi$ be the $d$-dimensional nontrivial irreducible representation of
$U_q(g)$ introduced earlier.  A $U_q(g)$-module action
\begin{align*}
   \Uq\otimes H&{\overset{\circ}{\longrightarrow}} H
\end{align*}
can be defined by setting
\begin{align}
   a\circ s_i = \sum_{j=1}^d \pi(a)_{ji} s_j.
\end{align}
Since $\pi$ is irreducible, $H$ contains a unique ( up to scalar multiples)
highest weight vector $s^\pi$.

Now each $H^r$ furnishes a $\Uq$-module with the module action defined by the
co-multiplication:
\[
   a\circ(s_{i_1} s_{i_2} \dots s_{i_r}) = \sum_{\{j\}} \pi(a_{(1)})_{j_1i_1}
      \pi(a_{(2)})_{j_2i_2} \dots \pi(a_{(r)})_{j_ri_r}
      \times s_{j_1}s_{j_2}\dots s_{j_r}.
\]
Now $\cal{H}$ becomes a $\Uq$-module, provided  that we consider only finite
linear combinations of vectors in the direct sum.

The dual vector spaces $H^{* r}$ have a  natural $U_q(g)$-module structure.
On ${\cal H}^*$, the $\Uq$ action is defined  by requiring
\[
   \langle a\circ s^*_i, s_j \rangle = \langle s_i^*, \gamma(a)\circ s_j\rangle
      = \pi(\gamma(a))_{ij}.
\]
Then each $H^{*r}$ becomes a tensor product module, and $\cal{H}^*$
 is the module obtained as the algebraic direct sum of the $H^{*r}$'s.
Explicitly the action of $\Uq$
on $H^{*r}$ is given by
$$
\langle a\circ(s^*_{j_r} s^*_{j_{r-1}} \dots s^*_{j_1}), s_{i_1} s_{i_2}
     \dots s_{i_r} \rangle
   = (\pi_{j_1i_1}\otimes \dots \otimes \pi_{j_ri_r}) \Delta^{(r-1)}
      (\gamma(a)).
$$
Note that if we write
\[
   \Delta^{(k-1)}(a) = \sum a_{(1)}\otimes ...\otimes  a_{(k)},
\]
then
\begin{align*}
   \Delta^{(k-1)}(\gamma(a)) &= \sum_{(a)} \gamma(a_{(k)}) \otimes
     \gamma(a_{(k-1)})\otimes ... \otimes \gamma(a_{(1)})
\end{align*}
Hence
\[
   a\circ\left(s^*_{j_r} s^*_{j_{r-1}}\dots s^*_{j_1}\right)
      = \sum_{(a)} (a_{(1)}\circ s^*_{j_r})(a_{(2)}\circ s^*_{j_{r-1}}) \dots
      (a_{(r)} \circ s^*_{j_1})
\]

\subsection{ $\Uq$ actions on
 $\cal{H}\otimes\cal{H}^*$ and $\cal{H}^*\otimes\cal{H}$}

The actions are defined by the co-multiplication in the obvious way, namely,
for $u\in{\cal H},\ v^*\in\cal{H}^*$,
\begin{align*}
   a\circ(u\otimes v^*) &= \sum a_{(1)}\circ u\otimes a_{(2)} \circ v^*, \\
   a\circ(v^*\otimes u) &= \sum a_{(1)}\circ v^*\otimes a_{(2)}\circ u.
\end{align*}
They have the following useful properties
\begin{align*}
    \sum_i a\circ (u s_i\otimes s^*_i v^*)
&\qquad = \sum_j \sum_{(a)} (a_{(1)}\circ u)s_j\otimes s^*_j(a_{(2)}\circ v^*), \\
 a\circ (v^* s^*_i\otimes s_j u)
   &\qquad = \sum_{k l} \sum_{(a)} (a_{(1)}\circ v^*) \pi(\gamma(a_{(2)}))_{i k }
   s^*_k\otimes \pi(a_{(3)})_{l j} s_l (a_{(4)}\circ u).
\end{align*}

We wish to examine the properties of the module actions under Cuntz multiplication.
Consider
\begin{align*}
   & a\circ\left(s^*_{j_r} s^*_{j_{r-1}} \dots s^*_{j_1} \otimes s_{i_1} s_{i_2}
      \dots s_{i_t}\right) \\
   &\qquad = \sum_{(a)} (a_{(1)}\circ s^*_{j_r})(a_{(2)}\circ s^*_{j_{r-1}})\dots
      (a_{(r)}\circ s_{j_1}^*)\otimes(a_{(r+1)}\circ s_{i_1}) \dots
      (a_{(r+t)}\circ s_{i_t}).
\end{align*}
Direct calculations can establish that
\begin{align*}
& a\circ\left(s^*_{j_r} s^*_{j_{r-1}} \dots s^*_{j_1}  s_{i_1} s_{i_2}
      \dots s_{i_t}\right) \\
   &\qquad = \sum_{(a)} (a_{(1)}\circ s^*_{j_r})(a_{(2)}\circ s^*_{j_{r-1}})\dots
      (a_{(r)}\circ s_{j_1}^*)(a_{(r+1)}\circ s_{i_1}) \dots
      (a_{(r+t)}\circ s_{i_t}),
\end{align*}
where both sides of the above equation are regarded as elements of the
Cuntz algebra.   It is also clearly true that
\begin{align*}
   & a\circ\left(s_{i_1} s_{i_2}
      \dots s_{i_t} s^*_{j_r} s^*_{j_{r-1}} \dots s^*_{j_1} \right) \\
   &\qquad = \sum_{(a)} (a_{(1)}\circ s_{i_1})(a_{(2)}\circ s_{i_2})\dots
      (a_{(t)}\circ s_{i_t})(a_{(t+1)}\circ s^*_{j_r})
     (a_{(t+2)}\circ s^*_{j_{r-1}}) \dots (a_{(r+t)}\circ s^*_{j_1}).
\end{align*}
Therefore, the $\Uq$ action preserves the Cuntz multiplication.

\subsection{$O_d^0$ as a module algebra over $U_q(g)$}

The results of the last subsection suggest that a $\Uq$-module action
on $O_d^0$:
\[
\Uq \otimes O_d^0{\overset{\circ}{\longrightarrow}} O_d^0,
\]
can be introduced directly in which each element of $\Uq$ acts by
a non-unital endomorphism (by which we mean it preserves the
multiplication but not the *-operation or identity) of
$O_d^0$. This is achieved by defining
\begin{align}
   a\circ 1& = \epsilon(a) \nonumber \\
   a\circ s_i&=\sum_{j} \pi(a)_{j i} s_j,\nonumber \\
   a\circ s^*_i&=\sum_{j} \pi(\gamma(a))_{i j}s^*_j,\nonumber \\
   a\circ (x y)&=\sum_{(a)}(a_{(1)}\circ x)(a_{(2)}\circ y),\quad x, y\in O_d^0.
\label{action}
\end{align}
The defining relations of the Cuntz algebra are clearly preserved,
\begin{align*}
   a\circ(s^*_i s_j) &= \sum_{(a)} (a_{(1)}\circ s^*_i)
                         (a_{(2)}\circ s_j)= \delta_{ij} a\circ 1,  \\
   a\circ (\sum_1^d s_i s^*_i)& = \sum_1^d \sum_{(a)} (a_{(1)}\circ s_i)
                             (a_{(2)}\circ s^*_i)=a\circ 1.
\end{align*}
Therefore, $O_d^0$ defines a module algebra over 
$U_q(g)$ under the action (\ref{action}).

This module algebra structure of $O_d^0$ over $U_q(g)$  can be 
straightforwardly extended to ${\tilde U}_q (g)$ by specifying the action 
of $k_\pi$ on the highest weight vector $s^\pi$ of $H$, 
\begin{align}
 k_\pi. s^\pi & = q^{(\mu,\ \mu)} s^\pi. \label{extension}
\end{align} 
Note that when the highest weight $\mu$ of $\pi$ does not 
belong to the root lattice,  $(\mu,\ \mu)$ is in general a rational 
number, and in that case we make a choice of the complex  value of 
$q^{(\mu,\ \mu)}$.

\subsection{Braids as fixed points}

In the previous subsection we established the first claim of Theorem 1.
To verify the remaining claims we
now consider the fixed points.
The fixed point set of $O_d^0$ under the $\Uq$  action (\ref{action})
is
\[ \{ u\in O_d^0  \quad  | \quad      a\circ u = \epsilon(a) u,
\quad \forall a\in\Uq\}. \]
Since $O_d^0$ is a $U_q(g)$ module algebra, the fixed point set
defines a subalgebra of $O_d^0$.
Consider the $\Uq$ action on the braid generator
$\theta = \sum \sigma_{ij,kl}  s_i s_k s^*_l s^*_j$
\begin{align*}
   a\circ \theta & = \sum \left[(\pi\otimes\pi)\Delta (a_{(1)})\cdot\sigma\cdot
      (\pi\otimes\pi)\Delta (\gamma(a_{(2)}))\right]_{ij,kl}
      s_is_k s^*_l s^*_j \\
   &= \sum\left[(\pi\otimes\pi)\Delta(a_{(1)}\gamma(a_{(2)}))\cdot\sigma
      \right]_{ij,kl} s_is_k s^*_l s^*_j\\
   & = \epsilon(a)\theta,
\end{align*}
where we have used the fact that the braid generator $\sigma$ commutes
with $(\pi\otimes\pi)\Delta(a)$, $ \forall a\in \Uq$.  More generally,
\begin{align*}
    a\circ \theta_{r+1}
   & = \sum a_{(1)}\circ(s_{i_1}\dots s_{i_r})\cdot [ a_{(2)}\circ
      \eta(\sigma) ]  a_{(3)}\circ(s^*_{i_r}\dots s^*_{i_1}) \\
   & =\epsilon(a)\theta_{r+1}.
\end{align*}
Thus we have shown that the braids  $\theta_r\in O_d^0$,
$\forall r=0,1,2,\dots,$
are fixed points of the $\Uq$-action.

\section{ $O_d^0$ as a co-module algebra over $\Gq$}

\subsection{The algebra $O_d^0$ as a co-module algebra over $G_q^{(\pi)}$}
Following the general method discussed in Section 2
we can define a co-module action of $G_q^{(\pi)}$ on $O_d^0$ by
\begin{align}
 \omega(1)&=1\otimes \epsilon,\nonumber     \\
 \omega(s_i)&=\sum_{j=1}^d s_j\otimes u_{ji},\nonumber     \\
\omega(s_i^*)&=\sum_{j=1}^d s_j^*\otimes \gamma_0(u_{ij}),\nonumber\\
%\qquad\qquad\quad
 \omega(x y)&= \omega(x) \omega(y), \qquad\qquad \forall \quad x , y \in O_d^0,
\label{coaction}
\end{align}
where the multiplication on the right hand side of the last equation is the
natural one for the algebra $O_d^0\otimes G_q^{(\pi)}$ induced by the
multiplication of $O_d^0$ and that of $G_q^{(\pi)}$.  The consistency of this
definition is confirmed by the simple calculations below:
\begin{align*}
\omega(x) \omega(y)(a)&=\sum_{(x),(y)}x_{(1)}y_{(1)}
           \langle x_{(2)} y_{(2)},  a\rangle \\
&=\sum_{(qa)}\sum_{(x),(y)}x_{(1)}\langle x_{(2)}, a_{(1)}\rangle
              y_{(1)}\langle y_{(2)}, a_{(2)}\rangle\\
&=\sum_{(a)} ( a_{(1)}\circ x) \  (a_{(2)}\circ y )=
%   \sum_{(a)}a_{(1)}\circ x\ a_{(2)}\circ y\\
a\circ(x y)\\
&=\omega(x y)(a),\qquad \forall a\in A.
\end{align*}
Hence $O_d^0$ is a co-module algebra over $G_q^{(\pi)}$.

\subsection{The braids as co-fixed points}
We have already shown that the braids are fixed by the $U_q(g)$ action, thus
they must be co-fixed by $G_q^{(\pi)}$.
Nevertheless, we look at an example to illustrate the general result.

Let
\[
\theta=\sum\sigma_{ij,kl}s_i s_k s_l^* s_j^*
\]
be a braid embedded in $O_d^0$.
The coefficient matrix  $\sigma$ commutes with $(\pi\otimes \pi)\Delta(a)$
$\forall a\in \Uq$.  This translates, for the dual Hopf algebra $G_q^{(\pi)}$,
to the following relations
\[
\sigma_{12}U_1U_2=U_1U_2\sigma_{12}.
\]
Now
\begin{align*}
\omega(\theta)
&=\sum s_is_k s_l^*s_j^*\otimes(U_1U_2\sigma_{12}\gamma_0(U_1U_2))_{ij,kl}\\
&=\theta\otimes\epsilon.
\end{align*}

The discussion of this section completes the proof of Theorem 1.

\subsection{Proof of Theorem 2}

As shown in 4.1, the linear span $H$ of the generators $s_i$ of
$O_d$ furnishes an irreducible $U_q(g)$ module.
It follows the results of 2.2 that there exists a nonvanishing
rank $d$
$q$-antisymmetric tensor $S_q\in H^d$, which generates a
trivial $U_q(g)$ module,
\[ a\circ S_q=\epsilon(a) S_q. \]
Clearly $S_q$ is also co-fixed by $\Gq$.
This $S_q$ is the distinguished element of Theorem 2.

To complete the proof of Theorem 2 we note that any element $x$ of $O_d^0$
which is fixed under the $\Uq$ action is also co-fixed under the
$\Gq$ action. For $0<q<1$ this means $x$ is co-fixed under the co-action of
$\overline\Gq$. Hence by \cite{P} $x$ is in the subalgebra generated
topologically by the
rank $d$ $q$-antisymmetric tensor together with the braid group elements.
On the other hand $x\in O_d^0$ and so lies in the algebraic subalgebra
generated by  the
rank $d$ $q$-antisymmetric tensor together with the braid group elements.

\section{ The $\Gq$ co-module $O_d^0$ as a $U_q'$ module }

In the above calculations we started with $O_d^0$ as a $\Uq$
module and showed how to obtain the dual co-action of $\Gq$.
It is of some interest to understand whether, when the reverse procedure is
adopted, namely regarding $O_d^0$ as a  $\Gq$ co-module as in \cite{P}
and constructing the dual action, one recovers the given $\Uq$ action.
This is indeed the case but it requires us to understand $\Uq$
in a different way namely in terms of $U_q^\prime$.

\subsection{Co-module of a Hopf algebra as a module of the dual Hopf algebra}

As in Section 2 we adopt a general viewpoint first
and specialise later.
In the notation of Section 2,
 $A$ is a Hopf algebra with the finite dual $A^0$,
which is also a Hopf algebra.  Let $V$ be a co-module of $A$,
\[
V{\overset{\omega}{\longrightarrow}} V\otimes A.
\]
The defining property of $\omega$ is that
\begin{align*}
(\id_V\otimes \Delta)\omega=(\omega\otimes \id_A)w:&
\quad V\rightarrow V\otimes A\otimes A, \\
(\id_V\otimes \epsilon)\omega = \id_V:& \quad V\rightarrow V.
\end{align*}
Using Sweedler's sigma notation  \cite{S}, we  write, for $v\in V$,
\[
w(v)=\sum_{(v)}v_{(1)}\otimes v_{(2)},
\qquad v_{(1)}\in V, \quad v_{(2)}\in A,
\]
where the right hand side is assumed to be a finite sum.

We observe that $V$ has a natural $A^0$-module structure,
\[
A^0\otimes V{\overset{0}{\longrightarrow}}V
\]
defined, for any $v\in V$ by,
\[
a\circ v=(\id_V\otimes a)w(v),\qquad a\in A^0
\]
or more explicitly,
\[
a\circ v=\sum_{(v)}v_{(1)}\langle a, v_{(2)}\rangle.
\]
To check that this indeed defines a module over $A^0$, consider
\begin{align*}
b\circ(a\circ v)
&=\sum_{(v)}v_{(1)}\langle b, v_{(2)}\rangle\langle a, c_{(3)}\rangle\\
&=\sum_{(v)}v_{(1)}\langle ba, v_{(2)}\rangle =(ba)\circ v.
\end{align*}
The unit of $A^0$ is $\epsilon$,
\[
\epsilon\circ v=\sum_{(v)}v_{(1)}\langle \epsilon, v_{(2)}\rangle=v.
\]

\subsection{$O_d^0$ as a module over $U_q^\prime\subset(G_q^{(\pi)})^*$}

The co-action of $G_q^{(\pi)}$ on $O_d^0$ is defined by (\ref{coaction}).
It follows from the discussion of the last subsection that this
co-action dualizes an action of $U_q^\prime$ on the Cuntz algebra:
\begin{align*}
a\circ 1=\epsilon(a),&\\
a\circ s_i=\sum_j s_j \langle a, u_{ji}\rangle, &\qquad
a\circ s_i^*=\sum_j s_j^*\langle a, \gamma_0(u_{ij})\rangle, \\
a\circ(x y)=\sum a(w(x y))&
=\sum x_{(1)}y_{(1)}\langle a, x_{(2)}y_{(2)}\rangle
=\sum_{(a)}(a_{(1)}\circ x)(a_{(2)}\circ y).
\end{align*}
By recalling the defining relations (\ref{definition}) of $U_q^\prime$,
we can see that $Ker\phi$ annihilates $O_d^0$.  Hence the actions of
$U_q^\prime$ and ${\tilde U}_q(g)$ on $O_d^0$ coincide,
\[  a\circ x = \phi(a)\circ x,   \forall a\in U_q^\prime, x\in O_d^0, \]
where  the action of ${\tilde U}_q(g)$ on $O_d$ appearing on the right hand side
is that defined by (\ref{action}) and (\ref{extension}). 
Therefore, we have recovered
the ${\tilde U}_q(g)$ action and hence {\it a fortiori} the $\Uq$ action
on the Cuntz algebra from the $G_q^{(\pi)}$ co-action.

\section{Categorical interpretation}

There is a way to interpret our earlier constructions more
abstractly. This may be useful in order to relate this
 paper to other examples.
Let $A$ be a Hopf $*$-algebra generated by the elements  of
$(u_{ij})_{i, j=1}^d$.    We assume that  $u$ is unitary,
that is,  $uu^*=I=uu^*$.
Let $H$ be a Hilbert space of dimension $d< \infty$, which furnishes
a co-representation of $A$,  $\omega: H\rightarrow H \otimes A$,
$\omega(e_i)=\sum_j e_j\otimes u_{ji}$.
Tensor products of $H$ with itself
yield co-representations of $A$ in a natural way.

We have not specified the relations
satisfied by the matrix elements of $u$ at this stage.
In fact, it is our purpose here to show that the defining
relations can be recovered from the co-algebraic structure
and the co-fixed points of the $A$ co-action on
the braided tensor category
generated by $H$.

Consider the tensor category $F_d$,  whose objects are $H^{\otimes r}$,
and arrows are the linear mappings
$T:H^{\otimes r}\rightarrow H^{\otimes s}$.
As in Doplicher \cite{D},  we associate to $F_d$ the Cuntz algebra $O_d$.
Then there is a natural $A$ coaction on $O_d$ defined by
\begin{align*}
   \Gamma &: O_d\to O_d\otimes A \\
   \Gamma(s_i) &= \sum_{j=1}^d  s_j\otimes u_{ji}.
\end{align*}
This co-action, when  restricted to the UHF-algebra $M_d^\infty$ of $O_d$
through the embedding $\eta(e_{i_1j_1}\otimes\dots\otimes e_{i_kj_k}) =
s_{i_1}\dots s_{i_k} s_{j_k}^*\dots s_{j_1}^*$,  gives
\[
   \Gamma_{\operatorname{UHF}}(e_{i_1j_1}\otimes\dots\otimes e_{i_kj_k}) =
      \sum \begin{Sb} a_1\dots a_k \\ b_1\dots b_k \end{Sb}
      e_{a_1b_1}\otimes\dots\otimes e_{a_kb_k}\otimes u_{a_1i_1}\dots
      u_{a_ki_k} u^*_{b_kj_k}\ldots u^*_{b_1j_1}
\]
for all positive integers $k$.

Set $u^k = u\otimes \dots\otimes u$.  Then for any
$x\in M_d^k$, $\Phi(x) = \Gamma_{\operatorname{UHF}}(x)
= u^k(x\otimes 1)(u^k)^*$.
The fixed points of $\Phi$ are those $x$ such that
$u^k(x\otimes 1)(u^k)^* = (x\otimes 1) \Rightarrow u^k(x\otimes 1)=
(x\otimes 1)u^k$.

Let $R$ be the solution of quantum Yang-Baxter equation such that
$R : H\otimes H\to H\otimes H$.
Let $P$ be the flip. Then $\sigma = \operatorname{PR}$ satisfies the
braid relations.  As elements of $O_d$, $\theta=\eta(\sigma)
= \sum \sigma_{ijkl} \eta(e_{ij}\otimes e_{kl})
= \sum \sigma_{ijkl} s_i s_k s^*_l s^*_j$ and
\[
   \theta_i = \eta(b_i) =
\sum_{\{l\}} s_{l_1}s_{l_2}\dots s_{l_{i-1}}\otimes\sigma\otimes
      s_{l_{i-1}}^*\dots s_{l_2}^* s_{l_1}^*
\]
satisfy the braid group relations.

Form then a braided tensor category $(F_d,R)$ by replacing
the commutativity of $F_d$ with the braiding $b_i$.
Denote by $(O_d)^\Gamma$ the fixed point algebra of $(F_d,R)$ under
the coaction $\Gamma$:
$(O_d)^\Gamma = \{x : \Gamma(x) = x\otimes 1 \}$.
Now we require that the fixed point algebra $(O_d)^\Gamma$ is generated
by the braids $b_i$. Then it follows from
$u^2(\sigma\otimes 1)= (\sigma \otimes 1)u^2$ that
\[
R(u\otimes 1)(1\otimes u) = (1\otimes u)(u\otimes 1)R,
\]
which is precisely the defining relations of a quantum group in the
formulation of \cite{Faddeev}.
 Therefore, the algebraic structure of $A$ is recovered
from the co-fixed point algebra $(O_d)^\Gamma$
associated to the braided category $(F_d, R)$.

Conversely,  starting with our
Hilbert space $H$ of finite dimension $d$ as an
irreducible  $U$ module, where $U$ is a Hopf $*$-algebra,
we can consider the tensor category whose
objects are $H^{\otimes n}$ and the arrows are linear mappings
 $H^{\otimes n} \to H^{\otimes n}$.
Assume there exists a nondegenerate pairing between $A$ and $U$ which
relates $\Delta,\gamma,\epsilon$ of $A$ with the
multiplication, unit of $U$,   and relates the algebra
structure of $A$ with the comultiplication and co-unit of $U$.
As for the antipode, we have
   $$\langle \gamma(x), a \rangle = \langle x, \gamma(a) \rangle,
   x\in A,  a\in U.
$$

As before $O_d$ is associated to the tensor category $F_d$. Then there
exists a natural $U$-module action on $O_d$ given by
$a\circ 1 = \epsilon(a)1,\ a\circ s_i = \sum_j \pi(a)_{j i}s_j,
\ a\circ s_i^* = \sum_j \pi(\gamma(a))_{i j} s_j^*,
\ a\circ(x y) = \sum_{(l)}(a_{(1)}x)\circ
(a_{(2)}y)\ \forall x, y\in O_d$.
This module action gives rise to a co-module action of $A$ on $O_d$
through equation (\ref{comodule}) and  the pairing
\[
   \langle u_{ij}, a \rangle = u_{ij}(a) = \pi(a)_{ij}.
\]
We can continue in this vein, reinterpreting the
earlier explicit constructions from the more abstract viewpoint.
However we desist after making one final comment.
If we also assume that there exists a braid operator $\sigma$ acting on
$H\otimes H$, such that
$[\sigma, (\pi\otimes\pi)\Delta(a)]=0, \ \forall a\in U$
then $\theta = \eta(\sigma)\in O_d$  can be easily
seen to be co-fixed by $A$.  Furthermore,
\[  R u_1 u_2 = u_2 u_1 R, \]
where $R=P\sigma$ and $u=(u_{i j})_{i , j=1}^d$.

\end{document}